# Third Party Risk Modelling and Assessment for Safe UAV Path Planning in Metropolitan Environments

Bizhao Pang [a], Xinting Hu [b], Wei Dai [a,c], Kin Huat Low [a,*]

[a] *School of Mechanical and Aerospace Engineering, Nanyang Technological University, Singapore 639798, Singapore*
[b] *School of Air Traffic Management, Civil Aviation University of China, Tianjin 300300, China*
[c] *Air Traffic Management Research Institute, Nanyang Technological University, Singapore 639460, Singapore*

**Abstract**: Various applications of advanced air mobility (AAM) in urban environments facilitate our daily life and public services. As one of the key issues of realizing these applications autonomously, path planning problem has been studied with main objectives on minimizing travel distance, flight time and energy cost. However, AAM operations in metropolitan areas bring safety and society issues. Because most of AAM aircraft are unmanned aerial vehicles (UAVs) and they may fail to operate resulting in fatality risk, property damage risk and societal impacts (noise and privacy) to the public. To quantitatively assess these risks and mitigate them in planning phase, this paper proposes an integrated risk assessment model and develops a hybrid algorithm to solve the risk-based 3D path planning problem. The integrated risk assessment method considers probability and severity models of UAV impact ground people and vehicle. By introducing gravity model, the population density and traffic density are estimated in a finer scale, which enables more accurate risk assessment. The 3D risk-based path planning problem is first formulated as a special minimum cost flow problem. Then, a hybrid estimation of distribution algorithm (EDA) and risk-based A* (named as EDA-RA*) algorithm is proposed to solve the problem. To improve computational efficiency, k-means clustering method is incorporated into EDA-RA* to provide both global and local search heuristic information, which formed the EDA and fast risk-based A* algorithm we call EDA-FRA*. Case study results show that the risk assessment model can capture high risk areas and the generated risk map enables safe UAV path planning in urban complex environments. The statistical analysis is also conducted to test the transformational impact of risk assessment model and risk-based path planning method. Obtained results show that the proposed risk assessment model and risk-based method are effective for all types of urban patterns.

## 1. Introduction

Applications of advanced air mobility (NASA 2020) have been extensively seeing in urban areas for various cases such as traffic monitoring, aerial photography, delivery, etc. Projection also shows that drone operation in metropolitan areas will continue to rise (Narkus-Kramer 2017). To handle the large-scale UAV operations with different tasks, autonomous flying capability is crucial. As one of the key enabler of autonomous flying, path planning problems have been widely investigated with purposes of minimizing flight distance and operational cost (Ha et al. 2018), energy consumption (Wai and Prasetia 2019), or maximizing coverage rate for surveillance mission (Wu, Wu, and Hu 2020). However, third party risk issues are essential for UAV operating in metropolitan areas, as UAV may fail and fall due to loss of control or navigation (Pang, Ng, and Low 2020). Falling UAV may cause fatalities to people (Koh et al. 2018) and damages to properties (Dalamagkidis, Valavanis, and Piegl 2008). UAV operating in low-altitude airspace also brings societal issues like noise impact and privacy concerns to the public (Lin Tan et al. 2021). These issues are considered as psychological risk cost of UAV operation to the public, which needs to be mitigated in path planning phase. In this paper, we investigate and answer the question of how to quantitively assess various UAV operational risks, and how to effectively mitigate these risks by using risk-aware airspace modelling and risk-based path planning method.

There are existing studies investigated the risk assessment problems of UAV operation, and they focused on impact probability and severity models to people and vehicle on the ground. The authors (Mitici and Blom 2019) presented main mathematic models for conflict and collision probability estimation, which provide insights for collision risk assessment of AAM. Pioneer works studied the probability of fatalities and the fatality rates associated with a ground impact to pedestrians, and analysis results showed that the risk of fatality to human is low in condition of light UAV operates in areas with low population density (Dalamagkidis, Valavanis, and Piegl 2008). The probability model of UAV to road traffic was also established. The authors (Bertrand, Raballand, and Viguier 2018) defined the possible ground impact area of falling UAV and developed the collision probability model, which helps for identification of main risky areas of road network. Follow up works studied the impact severity of UAV to people and they subsequently proposed the weight threshold of falling UAV impact ground people based on the injury scale and criterion (Koh et al. 2018; Clothier, Williams, and Hayhurst 2018). Based on the UAV impact probability and severity studies, researcher proposed a risk-based approach for small UAV operations (Breunig et al. 2019), and generated the probabilistic map using Monte Carlo simulation for more accurate ground impact risk analysis (Levasseur et al. 2019).

---

[*] Corresponding author: K.H. Low (mkhlow@ntu.edu.sg)

Recent studies paid attention to third-party risk modelling and analysis. The third-party risk was defined as risks pertaining to human life and property damage which are not onboard the UAV (Melnyk et al. 2014; Clothier, Williams, and Hayhurst 2018; Jiang, Blom, and Sharpanskykh 2020). In subsequent studies, a third-party risk framework was proposed to analyze the UAV ground impact risk (Melnyk et al. 2014), and third-party risk indicators and their utilization in safety regulations were proposed (Jiang, Blom, and Sharpanskykh 2020). By using these proposed frameworks, some practical studies have been conducted to model the third-party risk in urban environments (S. H. Kim 2020; Ren and Cheng 2020), and the level of risk for UAV system was also proposed to identify critical areas and actions (Gonçalves, Sobral, and Ferreira 2017). The analysis and modelling of these risks facilitate the generation of risk aware map, which can be used for risk-based path planning with aims of achieving safer UAV operation in metropolitan environments.

UAV path planning problems have also been extensively studied with different models and optimization objectives. Exact methods like Dijkstra algorithm (Dijkstra 1959), heuristic algorithms like A* (Bell 2009), and swarm-based heuristic methods (Wu 2021; Wu et al. 2021). These methods have various optimization objectives like minimizing travel distance and cost, maximizing flight duration. These methods always consider obstacle avoidance but ignore risks underneath the UAV flying path. Extended from conventional distance or cost based path planning problems, the risk-based one is relied on a risk map used for path planning (De Filippis, Guglieri, and Quagliotti 2011; Hu et al. 2020). Various methods and algorithms were developed to generate the risk map and to cope with the risk-based path planning problems. The A*-based algorithms, for instance, was developed together with Dubins Curves for risk-based path planning and smoothing (De Filippis, Guglieri, and Quagliotti 2011). In the follow up study, the authors (Primatesta, Guglieri, and Rizzo 2019) developed a RiskA* algorithm to minimize the risk of the produced path. Genetic algorithm and Dijkstra methods are also popular in addressing this problem and been chosen as benchmarks to compare with A*-based path planning algorithm (Da Silva Arantes et al. 2017; Votion and Cao 2019). Other methods like Markov decision process was used with hierarchical method to maximize efficiency and minimize risks (Feyzabadi and Carpin 2014). A rapidly exploring random tree (RRT) was proposed to minimize the third-party risk of UAV takeoff trajectories (Rudnick-Cohen, Azarm, and Herrmann 2019). And Tabu search algorithm was employed to optimize the UAV route in order to minimize the cost of damaged cargos (Zhu et al. 2020). What is more, authors (Chung et al. 2019) also developed a risk-aware graph search algorithm to select paths which have high probability to yield low risk. On the other hand, the UAV operational environments have also been covered from factory-like area (Feyzabadi and Carpin 2014) to inhabited areas (Rudnick-Cohen, Herrmann, and Azarm 2016) and to urban environments (N. Kim and Yoon 2019). In different areas, the risk types are various. For instance, in factory area the main risk sources are critical infrastructure and property damages. While in inhabited area and urban areas, population and vehicle density, high-rise buildings are more sensitive for risk-based path planning.

In overall, existing studies have investigated the probability and severity models of ground impact on human life and property. However, these models rarely considered the mobility of population density in metropolitan areas, which fails to accurately capture the population density distribution a major contributor in risk modelling. Various risk types (fatality, property, etc.) have been investigated individually in different environments. However, an integrated risk assessment model is still lacking to cope with various risks in complex metropolitan environments. Lastly, existing path planning methods rarely incorporated risk cost into fitness function, and with even less studies investigated the risk-based heuristic function to improve optimality and efficiency of the risk-based path planning method.

In this paper, we propose an integrated risk assessment model with 3D risk aware airspace modelling, and we also develop a robust and effective algorithm to address the risk-based path planning problems with the goal of minimizing operational risk. We summarize the main contributions of this article as follows.

(1) We establish an integrated risk assessment model with a gravity model to better estimate the population density distribution and to capture high risk areas in a finer scale. The model considers three main risk categories in urban environments, which includes fatality risk (human life), property damage risk (infrastructure), and societal impact risk (noise and privacy). The introduced societal impact risk enables public perception of drone operation been considered in safe and sustainable airspace management and UAV operation planning.

(2) We formulate the risk-based 3D path planning problem as a special case of minimum cost flow problem. The objective is to find a minimum cost flow starting from origin to destination (OD) among the graph, with constrains of motion step size, flight consistency and obstacle clearance.

(3) We develop a hybrid algorithm integrating estimation of distribution algorithm, k-means method and improved A* algorithm (named as EDA-FRA*) to solve large scale 3D risk-based path planning problems. The outer loop of the EDA-FRA* algorithm is a 0-1 optimization problem, which aims for selecting and optimizing the low risk-cost path points based on OD information. The k-means clustering algorithm is introduced to extract heuristic information from selected low-risk path points for A* path searching algorithm to produce a risk-cost-effective path with high robustness and efficiency.

The rest of the paper is structured as follows. Section 2 analyzes the risk types in urban environments and illustrates the concept of risk aware airspace and path planning. The integrated risk assessment model is established in Section 3. The mathematic problem formulation and the hybrid 3D path planning algorithms are developed in Section 4. Followed by simulation validations and case studies in Section 5. Section 6 concludes the main findings of this article.

## 2. Problem background

In metropolitan environments, there are dense populations, high-rise buildings, critical infrastructure, etc. UAV operates in such low altitude airspace will encounter various risk issues (Ghasri and Maghrebi 2021). In this paper, the scope of operation altitude is below 400 feet (Federal Aviation Administration 2020) above the ground. Recent studies investigated various risks in urban environments, and most of their attentions are on the risks of the impact on ground people and road network (Bertrand, Raballand, and Viguier 2018; Clothier, Williams, and Hayhurst 2018), midair collision with small UAS and manned aircraft (Zou et al. 2021; Wang, Tan, and Low 2019), impact of noise and privacy issues (Vascik and Hansman 2018; Lin Tan et al. 2021), as well as UAV operational cost and efficiency (Ha et al. 2018). We conclude the primary risk sources as three categories as follows, illustrated in Fig. 1.

(1) Fatality risk. UAV impact pedestrians and vehicles on the ground, causing injuries or fatalities on people.
(2) Property damage risk. Falling UAV hits critical infrastructures or collides with high-rise buildings, causing property loss.
(3) Societal impact risk. Noise and privacy impact to the public is a big concern for the acceptance of UAV operation in urban environments. These impacts are modelled as societal impact risk, which will be assessed and mitigated.

Other risk factor like midair collision of UAV and manned aircraft in integrated urban airspace (Vascik and Hansman 2019) is not considered in this work, because airport performs segregated operation with UAV and the Aerodrome Control Zone is treated as restricted area where the UAV is strictly not allowed to enter. Risks of UAV intruding military-related bases and facilities are also out of scope of this article.

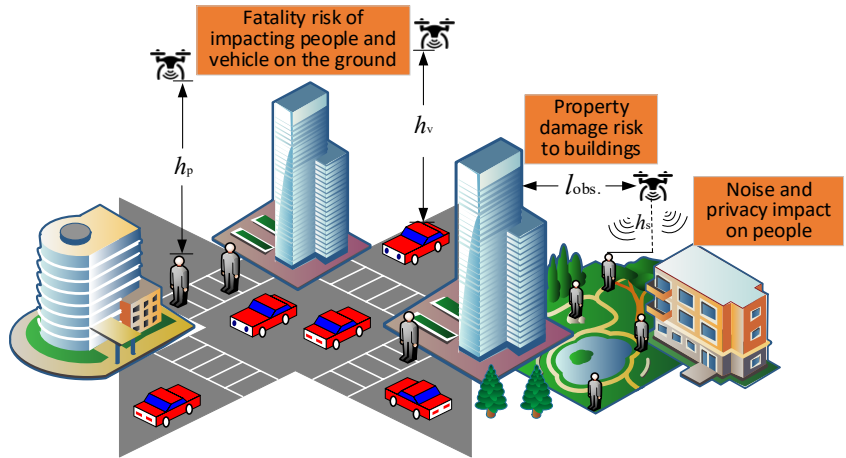

**Fig. 1.** Illustration of three primary risk types in metropolitan environments.

The identified risk sources of population density, vehicle density and buildings are discretely distributed in metropolitan environments. Subdivision of airspace into smaller manageable unit enable more flexible managing (Cho and Yoon 2019). To quantitatively assess these risks, urban low-altitude airspace is divided into discrete 3D air block unit (Pang, et al. 2020) and the centroid of the unit is denoted as $v_{xyz}$ (Fig. 2(a)). UAV operates from one air block to another with 26 possible points to choose for the next move. Risk assessment of each airspace unit is conducted based on its pertaining environments such as population density and vehicle density underneath. UAV operates following the centroid point to maintain safe separation with other UAVs in adjacent air blocks. The risk is represented as colored air blocks (Fig. 2(b)), and the 3D risk map is illustrated as Fig. 2(c). UAV operates in complex 3D risk map to avoid high risk areas (presented as red color) and to minimize total operational risk.

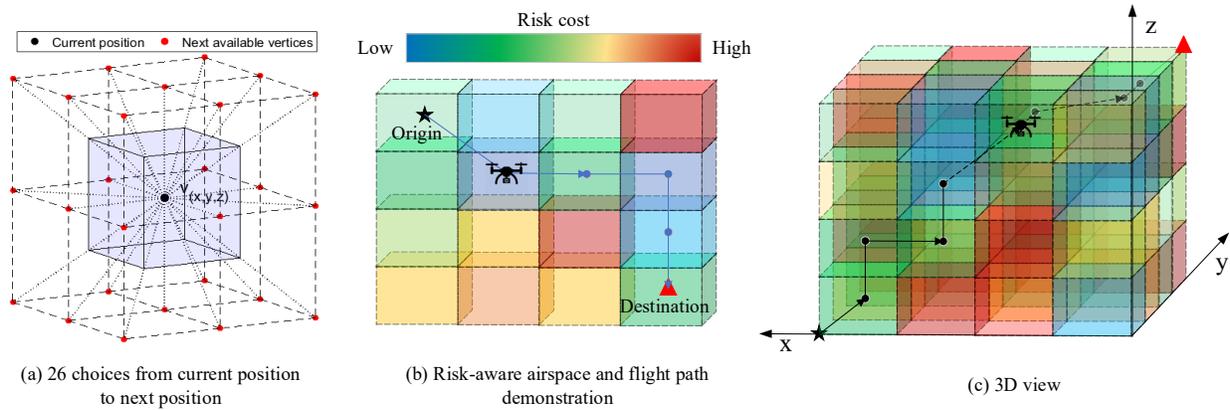

(a) 26 choices from current position to next position

(b) Risk-aware airspace and flight path demonstration

(c) 3D view

**Fig. 2.** Risk aware airspace modelling and risk representation.

The core of the risk aware airspace model and path planning is the accurate and quantitative risk value of each airspace unit. To achieve this, an integrated risk assessment model is proposed to compute the risk value.

## 3. Risk modelling and assessment

The integrated risk assessment model includes three main parts: fatality risk cost model, property damage risk cost model, and societal impact risk cost model.

### 3.1. Fatality risk cost model

#### 3.1.1. Risk of UAV impacts people on the ground

As it is possible that UAV might be loss of control or power when operating, falling UAV may impact people on the ground (see Fig. 1). There are three processes (Bertrand et al. 2017; Hu et al. 2020) a crash incident will cause injury or fatality to pedestrians: (a) failure of UAV; (b) falling UAV impacts people on ground; and (c) fatality damage caused to the people.

That falling UAV hits people on the ground causing fatalities is a chain action, which corresponds with the three processes mentioned above. The risk cost of UAV impacting ground people is defined as the number of fatalities per hour, denoted as

$$C_{r\_p} = P_{crash} N_{hit}^{p} R_{f}^{p} \qquad (1)$$

where $C_{r\_p}$ is the risk cost associated to the fatality of people, and $P_{crash}$ is the probability of UAV failure. Note that $N_{hit}^{p}$ is the number of pedestrians hit by falling UAV (proportional to the population density of people), and $R_{f}^{p}$ is the fatality rate associated to the function of kinetic energy.

The $P_{crash}$ is primarily determined by the capability of UAV itself, including hardware and software capabilities and reliability. The $R_{f}^{p}$ is strongly correlated with the weight and falling height of the UAV. The most uncertain variable in Eq. (1) is the $N_{hit}^{p}$, which associates with the population density, defined as

$$N_{hit}^{p} = S_{hit} \sigma_{p} \qquad (2)$$

where $S_{hit}$ is the explored area of UAV impacts the ground, and $\sigma_{p}$ is the population density in the administrative unit $u$.

The fatality rate $R_{f}^{p}$ associates with two main factors: impact kinetic energy and sheltering effects. The kinetic energy $E_{imp}$ of falling UAV primarily determines the severity of impact, while the sheltering coefficient $S_{c}$ affects the degree of impact on the people and vehicle, as the buffering effects of buildings, trees, etc. will soften the ground impact on them. Inspired by (Primatesta, Rizzo, and la Cour-Harbo 2020), the sheltering coefficient $S_{c}$ is introduced as the absolute real number $S_{c} = (0,1]$, and the fatality rate is presented as

$$R_f^p = \frac{1}{1+\sqrt{\frac{\alpha}{\beta}}(\frac{\beta}{E_{\text{imp}}})^{\frac{1}{4S_c}}} \tag{3}$$

where $\alpha$ is the impact energy that might cause 50% fatality with $S_c = 0.5$, while $\beta$ is the impact energy threshold required to cause fatality as $S_c$ approaching zero (see Fig. 2 in (Dalamagkidis, Valavanis, and Piegl 2008)). Based on that, we take $\alpha = 10^6$ J and $\beta = 100$ J.

The impact kinetic energy $E_{\text{imp}}$ of the falling UAV is known as

$$E_{\text{imp}} = \frac{1}{2}mv^2 \tag{4}$$

where $m$ (kg) is the mass of the falling UAV, and $v$ is the velocity when UAV hitting the ground stuff. To compute $v$, we have the followings.

The vertical drag force $F_d$ of falling UAV is related to its size and materials, as well as the density of air, etc., denoted by (Koh et al. 2018)

$$F_d = \frac{1}{2}R_I S_{\text{hit}} \rho_A v_{\text{TAS}}^2 \tag{5}$$

where $R_I$ is the drag coefficient related to the UAV type, $\rho_A$ is the density of air (1.225 kg/m³ at sea level), and $v_{\text{TAS}}$ is the true air speed of falling UAV.

Then the acceleration of UAV is

$$a = \frac{F_g - F_d}{m} = g - \frac{R_I S_{\text{hit}} \rho_A v_{\text{TAS}}^2}{2m} \tag{6}$$

where $F_g$ is the gravitational force, $F_g = mg$. (g=9.8m/s²)

Thus, the $v$ at moment $t$ which UAV hits ground can be obtained as

$$v = \int_0^t \left(g - \frac{R_I S_{\text{hit}} \rho_A v_{\text{TAS}}^2}{2m}\right) dt = \sqrt{\frac{2mg}{R_I S_{\text{hit}} \rho_A}\left(1 - e^{-\frac{hR_I S_{\text{hit}} \rho_A}{m}}\right)} \tag{7}$$

where $h$ is the start falling height of UAV above the ground surface.

### 3.1.2. Risk of UAV impacts vehicle on the ground

Similar to the risk model of people impact, there are also three components of a crash incident on road network (Bertrand, Raballand, and Viguier 2018): (a) UAV failure; (b) falling UAV hits a ground vehicle; (c) the crash incident causes a traffic accident which subsequently cause injuries or fatalities to people.

The expected fatality of UAV impacting a ground vehicle can be defined as the number of fatalities per hour caused by falling UAVs, denoted as

$$C_{r\_v} = P_{\text{crash}} N_{\text{hit}}^v R_f^v \tag{8}$$

where $C_{r\_v}$ is the risk cost associated to vehicle. Note that $N_{\text{hit}}^v$ is the number of vehicles hit by falling UAV (proportional to the traffic density), and $R_f^v$ is the average fatality rate associated to vehicle accident.

The average number of ground vehicles which may hit by falling UAV can be defined as the ratio of total area of all vehicles projected and the total road area, denoted as

$$N_{\text{hit}}^v = S_{\text{hit}} \sigma_v \tag{9}$$

in which $\sigma_v$ is the traffic density in the administrative unit $u$.

### 3.1.3. Estimation of population density and traffic density using gravity model

The population density and traffic density distributions in metropolitan environments are the essential variables which will directly influence the UAV operational risk costs as discussed in Eq. (1), (2) and Eq. (8), (9). Based on the previous studies, these density distributions are strongly correlated with the consumption amenities (Rappaport 2008), which attracts people and vehicle. To quantitively assess this correlation between consumption amenity and population density, the gravity model was used to calculate the population density (Yao et al. 2017). Inspired by the gravity model (Pang et al. 2021) and population mapping method (Deville et al. 2014), we have the following formulas to compute the population density in urban environments.

The population density of given unit $u$ is given as

$$\sigma_p = e^{(1-r^2)} \sigma_{p.avg} \qquad (10)$$

where $\sigma_{p.avg}$ is the average population density in the given area. Note that $r$ is the radius of the gravity influence area induced by the amenity, which is given as 1 km in this work. As shown in Fig. 3, the population density decreases in an inverted exponential pattern with increase of the radius $r$. In first 0.3 km, the index remains high, capturing the high population density distribution in the very vicinity of amenities. While in range of 0.3 km to 1.0 km, the index drops linearly, demonstrating the even decrease of population density.

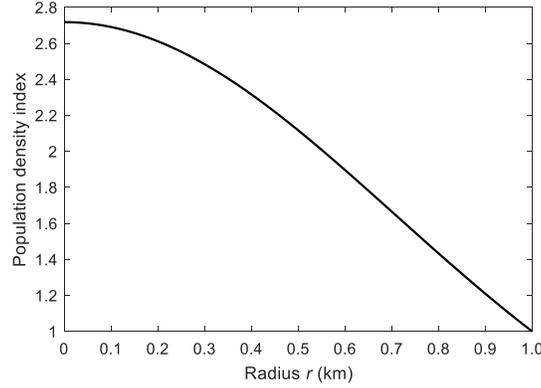

**Fig. 3**. Illustration of population density index changes with influence radius $r$.

Similarly, the road traffic density distribution can be denoted as

$$\sigma_v = e^{(1-r^2)} \sigma_{v.avg} \qquad (11)$$

where $\sigma_{v.avg}$ is the average traffic density in the given area.

The UAV operational risks to people and vehicle can be considered as fatality risk cost $C_{r\_f}$, presented as:

$$C_{r\_f} = C_{r\_p} + C_{r\_v} \qquad (12)$$

### 3.2. Property damage risk cost model

Dense high-rise buildings in urban environments is another challenge to perform UAV operations. Potential collisions with buildings pose property damage risks, and densely distributed high-rise buildings also limit the speed of traffic flow, resulting in the inefficiency of UAS system (Ang and Hansen 2019). Thus, the property damage risk cost model also integrates the operational efficiency cost, which are accounted for planning and optimization of airspace and traffic flow.

The flight altitude is a primary variable of the property damage risk model. As Fig. 4(a) depicted, in low altitude layer (Layer 1 for instance), the density of building is high. UAV operating in Layer-1-type airspace needs to frequently perform deconflictions to avoid obstacles, thus increasing risks and efficiency loss. In high altitude layer (i.e. Layer 4), in contrast, there are few buildings to affect UAV operation, so that the operational safety and efficiency can be significantly improved.

The building height distribution is not fit with standard normal distribution but log-normal distribution (Kirtner and Anderson 2008; Usui 2019), as building height is the nonnegative value and its distribution is not symmetry. Based on the height distribution relationship, the correlation between building height and property damage risk cost can be established as

$$f(h; \mu, \sigma) = \frac{1}{h\sigma\sqrt{2\pi}} e^{-\frac{(\ln h - \mu)^2}{2\sigma^2}} \qquad (13)$$

$$C_{r\_p.d} = \begin{cases} f(e^\mu), & 0 < h \leq e^\mu \\ f(h), & h > e^\mu \end{cases} \qquad (14)$$

where $\mu$ and $\sigma$ are the mean and standard deviation of the logarithmic variable (building height $h$). Note that $C_{r\_p.d}$ is the risk cost of property damage upon drone operation. For buildings with height smaller than the threshold of $e^\mu$, the risk cost equals to the one which the height $e^\mu$ has (as Eq. (14) defines). Which is because below that height ($h = e^\mu$), buildings are dense and the risks are high. The biggest risk cost value is therefore being given, and which is taken at height $h = e^{\mu}$. In this case, $\mu = 3.0467$ as computed above. While for buildings with height greater than $e^{\mu}$, the operational risk cost is computed as the log-normal distribution presents in Eq. (13). Meaning that with the increase of building height ($h > e^{\mu}$), the property damage risk cost decreases, as in higher layers there are few building obstacles to influence the safety of UAV operations. Note that the property damage risk factor is to facilitate the determination of optimal flight layer in particular areas. UAV should not collide with buildings in any flight layer with any building density.

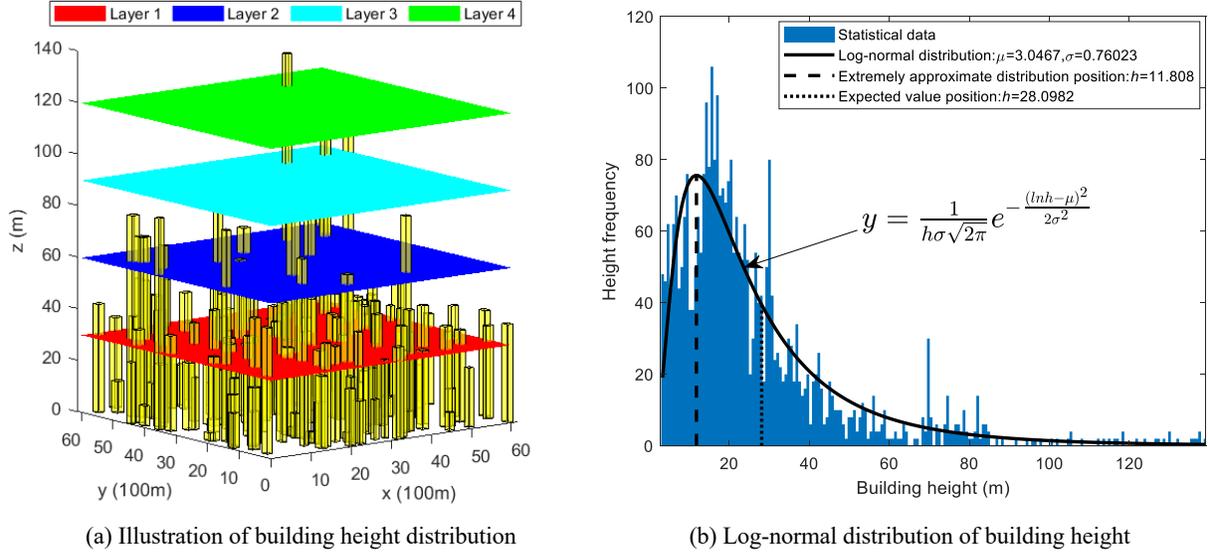

(a) Illustration of building height distribution  (b) Log-normal distribution of building height

**Fig. 4.** Building obstacles impact on UAV operation in urban area. In (a), the building height data is selected from a particular city area in Singapore, where 100 building clusters have been selected to illustrate the building height distribution. In (b), the statistical features of building height distribution are analyzed using log-normal distribution. Here, the height frequency presents the number of buildings at such height.

*3.3. Societal impact risk cost model*

Noise and privacy impacts are important societal issues and need to be considered when UAV operates in low altitude urban environments (Lin Tan et al. 2021), as low-flying UAV may upset people and make them feel annoying. That will be therefore considered as risk cost when conducting planning. The impacts of noise and privacy issues to the public are the same because their impacts are in effective when UAV operating close to people especially at nighttime. While with the increase of flying altitude, the impacts will decrease to the threshold which will not have effects on ground people. As the privacy risk cost is hard to be captured while it has the same nature of impact with noise issue, the societal risk model is therefore presented by noise impact risk model. The correlation of noise induced risk cost and its flying height is illustrated as Fig. 5. Based on the analysis, we know that the key factor of noise impact to people is UAV flying height.

A good first approximation of sound propagation is the spherical spreading, denoted as

$$I(si) = \frac{1}{h^2+d^2} \quad (15)$$

where $I(si)$ is the sound intensity at height $h$ and distance $d$ from the point directly under the drone. Here $d$ is taken as 30 feet (Alexander and Whelchel 2019).

$$L(sl) = \varpi L_h I(si) \quad (16)$$

in which $L(sl)$ is the sound level (dB); $\varpi$ is the convert coefficient from sound intensity to sound level; $L_h$ is the reference noise produced by drone, taken as $L_h = 60dB$ (Bulusu et al. 2017).

$$C_{r\_n} = L(sl) = \varpi L_h \frac{1}{h^2+d^2} \quad (17)$$

where $C_{r\_n}$ is the risk cost of noise upon drone operation in the given airspace unit. Noise impact will not be considered as risk cost for UAV operation if flying height exceeding the threshold. Based on previous studies (Bauer 2019; Torija, Li, and Self 2020), we take the height threshold as the one corresponded to the noise level of 40dB, illustrated as Fig. 5(b).

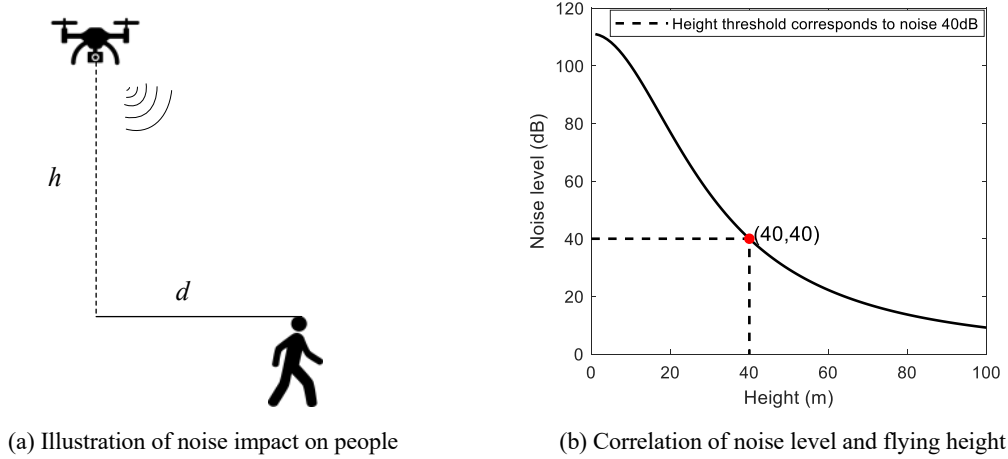

(a) Illustration of noise impact on people  (b) Correlation of noise level and flying height

Fig. 5. UAV noise impact on people.

*3.4. Integrated risk cost model*

The three risk cost models are discussed and developed above. To integrated them together as a comprehensive model, normalization is made for each type of risk. Obtained risk cost in each airspace unit $u$ will be divided by the maximum risk value of their own categories. All risk cost values of each type will be therefore ranged in (0,1]. The generalization factors are the reciprocal of the maximum risk cost for the corresponding type, denoted as:

$$\begin{cases} \omega_{r\_f} = \frac{1}{C_{r\_f}(max)} \\ \omega_{r\_p.d} = \frac{1}{C_{r\_p.d}(max)} \\ \omega_{r\_n} = \frac{1}{C_{r\_n}(max)} \end{cases}, \{C_{r\_f}(max), C_{r\_p.d}(max), C_{r\_n}(max)\} > 0 \qquad (18)$$

Note that $\omega_{r\_f}$, $\omega_{r\_p.d}$ and $\omega_{r\_n}$ are the generalization factors, which are used to keep the same order of magnitude for the three risk types.

The total operational risk cost in the given airspace unit integrates fatality cost, property damage risk cost and societal risk cost together. As the weight of these three types of cost might be different due to their significance or user's preferences (Liu et al. 2020), the contribution of each type of cost will also be various. For instance, aviation regulators may take safety as top priority, requiring a very low fatality risk cost of UAV operation. In this regard, the weight of fatality cost will be much greater than the other two factors. Thus, areas with dense population and vehicle will be identified as high-risk areas by proposed model, and the path planning will subsequently avoid these areas. To quantify the significance and preferences of UAS stakeholders on different risk types, the weight factors are introduced, and the total operational risk cost is computed as

$$C_{r\_total} = \sum_{i=1}^{3} \alpha_i \omega_i C_i \qquad (19)$$

where $\alpha_i$ is weight factor, $\omega_i$ is generalization factor and $C_i$ is risk cost. Here $i = 1, 2, 3$ corresponds to fatality risk, property damage risk and noise impact risk, respectively. Note that $\{\alpha_1, \alpha_2, \alpha_3\} \in [0,1]$, while $\alpha_1 + \alpha_2 + \alpha_3 = 100\%$.

## 4. Modelling of risk-based 3D path planning

In defined discrete airspace environment, the risk-based path planning problem is formulated, based on graph theory, as minimum cost flow problem in an undirected graph. To solve this problem, a hybrid algorithm has been proposed incorporating EDA and A* algorithms. The outer loop of the EDA-RA* algorithm is a 0-1 optimization problem, which aims for selecting and optimizing the path points as feasible search region. The optimized feasible region will be feed into the inner loop, where A* algorithm is employed, to generate the cost-effective path. To better improve the computational efficiency of the hybrid algorithm, the k-means clustering algorithm is introduced and incorporated to provide heuristic information for A* path searching algorithm, which is named as EDA-FRA*. Detailed problem formulation and algorithms development are presented in followings.

*4.1. Problem formulation of risk-based path planning*

To facilitate the risk cost assessment, the airspace is divided as uniform unit $u$ in three-dimensional space. The centroid of the airspace unit $u_i$ is presented as vertex $v_i(x_i, y_i, z_i)$, and the UAV operational risk cost in that unit is denoted as $C_{r\_i}$. The problem of finding a path from origin to destination with minimum cost is a special case of the minimum cost flow problem, which can be modelled as follows. Let $G = (V, E)$ be an undirected, connected and weighted graph such that all edge weights are nonnegative, with weight function $C: E \to \mathbb{R}_0^+$ and let $s$ and $t$ be distinct vertices of $G$. A path $P$ from $s$ to $t$ in $G$ is called the most risk-cost-effective path if $C(P) = \sum_{e \in P} C(e)$ is minimum among all paths from $s$ to $t$ in $G$. Here the weight function is equivalent to the risk cost function of Eq. (19).

A. Objective

The objective of this work is to minimize the total risk cost of planned path for UAV operation. The total risk cost consists of human fatality risk, property damage risk and societal impact risk.

$$\min: C(P) = \sum_{e_j \in P} C(e_j), \ j > 0 \tag{20}$$

where $C(P)$ is the total risk cost of the path $P$, and $C(e_j)$ is the risk cost of edge $e_j$. Note that P is the set of all edges included in the path.

Based on the risk assessment model, minimizing the total risk cost is to optimize several key variables, which are positions and flight altitude of UAV. In the model, they are presented as the 3D coordinate $(x_i, y_i, z_i)$, for each path point.

B. Constrains

As Fig. 2(a) shows, there are 26 available vertices can be chosen as the next point to move. Specifically, there are 6 vertices which are straightly connected with the current vertex $v_i$, while 12 vertices connected as planar diagonal and 8 vertices connected as cubical diagonal. Let $l$ be the length of the unit $u$. The motions and constrains of UAV in the discrete airspace can be expressed as follows.

$$\begin{cases} x_{i+1} = x_i + \xi_x^i, x_i \in \{0, l, 2l, \dots, X_{max}\} \\ y_{i+1} = y_i + \xi_y^i, y_i \in \{0, l, 2l, \dots, Y_{max}\} \\ z_{i+1} = z_i + \xi_z^i, z_i \in \{0, l, 2l, \dots, Z_{max}\} \\ \quad s.t. \begin{cases} \xi_x^i, \xi_y^i, \xi_z^i \in \{-l, 0, l\} \\ |\xi_x^i| + |\xi_y^i| + |\xi_z^i| \neq 0 \end{cases} \end{cases} \tag{21}$$

where $(x_{i+1}, y_{i+1}, z_{i+1})$ is the next path point $v_{i+1}$ we are choosing to move, and $(x_i, y_i, z_i)$ is the current position $v_i$. Note that $\xi_x^i, \xi_y^i, \xi_z^i$ are the unit lengths of each move corresponding to x-axis, y-axis and z-axis, respectively. Here the $X_{max}, Y_{max}$ and $Z_{max}$ are the boundaries of defined airspace in each axis.

As the step size in each axis is $l$, the first constrain is the unit moving length, which can only be chosen from alternative values of $\{-l, 0, l\}$, presented in Eq. (21). Assuming hovering is not allowed, the second constrain is that UAV must then take a move in whatever axis, which is denoted as the sum of the unit length must not equal to zero.

The third constrain related to obstacle avoidance. UAV should not collide with buildings in any situation, meaning that UAV should not enter the airspace unit, which is occupied by buildings. Therefore, we give the infinite risk cost to points which belong to building-occupied airspace, presented as

$$\{C_{r\_k} = \infty : u_k \in U_{Obs}, k \geq 0\} \tag{22}$$

in which $C_{r\_k}$ is the cost of these path points which are included in the occupied airspace. Note that $U_{Obs}$ is the airspace unit set containing all building-occupied units $u_k$.

*4.2. A hybrid EDA-RA\* algorithm for risk-based path planning*

To solve the developed minimum cost flow problem, there are several types of algorithms we can choose. Exact methods, like Dijkstra algorithm, is one of the classic and effective graph-based path-searching methods, which has been extensively used for path planning problems. However, Dijkstra is a computational inefficiency method, especially in dealing with large-scale problems like the one this paper studied. Heuristic methods, like A\* algorithm, has better performance in solving path planning problem in terms of computational time without reducing the quality of solutions, provided the accurate heuristic information can be offered.

For the standard A\* algorithm, the heuristic distance can be surely determined either as Manhattan distance or Euclidean distance (Fig. 6(a)) to estimate the distance from current position to destination. However, in the risk-based environment, the cost of each grid is different and unevenly distributed, making it hard to determine the heuristic distance (Fig. 6(b)). What is more, as the scale of the problem getting large, it is difficult for the heuristic methods which initiated with only one solution to search for the outstanding solutions from the feasible ones. In this regard, conventional A* algorithm is not suitable for the problem.

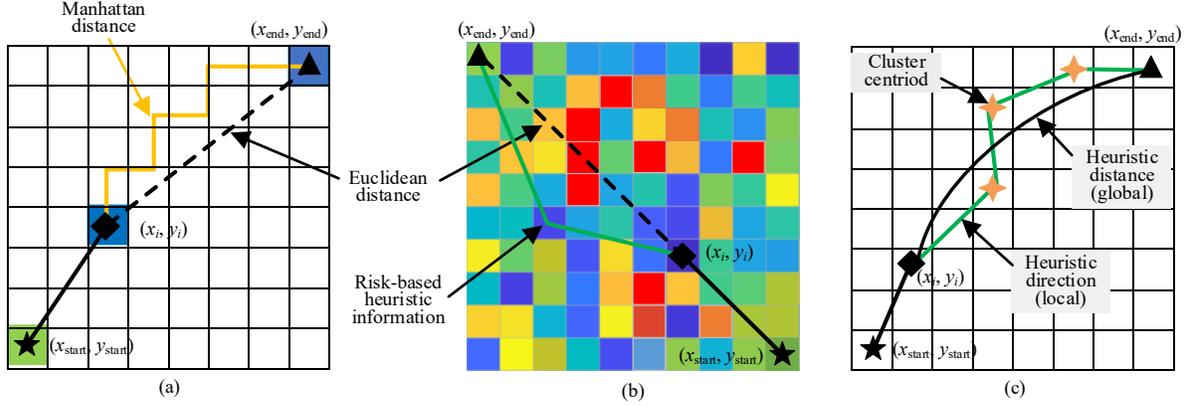

**Fig. 6.** Analysis of heuristic distance between classic A* algorithm and risk-based A*. Note that (a) shows that the classic A* algorithm normally has two options of its heuristic distance, Manhattan distance and Euclidean distance. For (b), the heuristic distance is chosen as risk-based distance, as the risk cost distribution from current node $(x_i, y_i)$ to target node $(x_{end}, y_{end})$ is unknown and unpredictable. For (c), the heuristic distance path is directly from current node to destination (global heuristic information). While the heuristic direction path is generated by connecting the cluster centroids to provide local heuristic information.

Swarm-based heuristic methods, on the other hand, initiate with a number of solutions and are suitable to solve the outer loop of our proposed problem. EDA is one of the typical swarm-based algorithms for solving both the continuous and discrete optimization problem. In this paper, the optimization variables of the problem are the 3D coordinate of path points, and the number of them are not fixed. EDA method is well aligned with these requirements, as it has no limitation for the number of variables, and it performs well in terms of global searching. In this work, the EDA method will be used to solve the outer loop 0-1 optimization problem, and it is incorporated with improved A* algorithm to generate the cost-effective path.

EDA algorithm is a stochastic method. The core of it is to generate and sample explicit probabilistic models of the promising solutions to guide the search for the optimum. The optimization process can be seen as a series incremental update of the probabilistic model to achieve the global optimal solution. That characteristic is good for globally optimizing the feasible region for the A* searching algorithm, which is the inner loop of the hybrid method to generate the path point by point. The general outline of the hybrid EDA-RA* algorithm for min-cost path planning is presented as **Algorithm 1**.

The EDA-RA* algorithm operates on a 3D operational risk-based airspace map, and the cost value is computed by the proposed risk cost assessment model. The output of this algorithm is the optimized path and total risk cost of the path. The main loop of the algorithm is from Line 6 to Line 23, which is concerned with the EDA method to optimize the feasible region, and A* algorithm for path generation.

The essential part of EDA-RA* algorithm is the probability update function (Line 22). By selecting the dominant populations from the species, the function is updated towards that the individuals, which belongs to dominant populations, will have increasing probabilities to be selected as optimal points in the graph. The update function is denoted as:

$$\boldsymbol{p}_{i+1} = (1 - l_{\text{rate}})\boldsymbol{p}_i + l_{\text{rate}} \frac{\boldsymbol{DS}}{DN} \tag{23}$$

where the $\boldsymbol{p}$ is the array that stores the probability of being selected for each individual, and $l_{\text{rate}}$ is the learning rate, which is the evolving factor to accelerate the optimization process with the accumulation of dominant species data. Note that $\boldsymbol{DS}$ is the array of dominant species, and DN is the total number of the species.

The selection of dominant populations is based on the fitness value. Here, the fitness value is the sum of risk cost value for vertices belong to the selected path, presented as

$$f_{\text{value}} = \sum_{v_i \in \text{path}} C_{r\_v_i} \tag{24}$$

As EDA is a stochastic method based on probabilistic model, during the iteration, some populations might be stuck by obstacles and may not be able to find path from origin to destination. The total risk cost for these populations will be replaced by the maximum risk cost among the whole species (Line 18). In subsequent iterations, these populations will be eliminated, and the selected dominant population will not have problem to find a feasible path.

**Algorithm 1**: Hybrid EDA-RA* algorithm for minimum risk-cost 3D path planning

1: Result: *path*, *TotalCost*
2: load CostDataset
3: *path* =*null*;
4: *TotalCost* =0;
5: initialize the probability matrix for EDA
6: **for** *i*=1:iterations %EDA outer loop
7:    **while** *j*<=populationSize
8:       r=rand(size(CostDataset));
9:       species{*j*,1} = 1.*(r<probability);
10:      *j*=*j*+1;
11:    **end**
12:    save species
13:    **for** k=1:populationSize % A* inner loop
14:       path=A*(species, obstacle);
15:       TotalCost=FitnessValue(path);
16:       TC=[TC;TotalCost];
17:    **end**
18: TC(replaceNP, :)=max(TC);
19: FitnessValue=TC;
20: [Fitness, index] = sort(FitnessValue);
21: dominantSpecies{:, 1} = species{index(:), 1}; %select dominant population
22: probability = (1-$l_{rate}$)*probability+$l_{rate}$*dominantSpecies/dominantNum;
23: **end**
24: **Return** (*path*, *TotalCost*)

*4.3. An improvement of EDA-RA\* with fast computation: EDA-FRA\**

The EDA-RA* algorithm incorporates A* as inner loop, and it is conducted for every population of the species at every iteration. Although the A* conducts very fast, this algorithm will cost considerable computational time as the problem scale getting larger. To cope with this problem, we further improve the EDA-RA* algorithm by introducing k-means method to provide both global and local heuristic information for path searching (see Fig. 6(c)). The improved hybrid algorithm is named as EDA-FRA*. Besides, the A* algorithm will also be improved to cater for the unique needs of risk-based 3D path planning problem.

The hybrid EDA-FRA* algorithm has three main functions. First function is EDA algorithm, which is used to globally optimize the feasible region that has low risk cost among all searching space. The second is the k-means algorithm. Based on the optimized feasible region, k-means clusters the feasible vertices to identify the heuristic directions (main tracks) and heuristic distance factor. The identified heuristic information will be ingested into the improved risk-based A* algorithm (named as RiskA* for easy reference) to generate the risk-cost-effective path. As the improved RiskA* algorithm is only called once, the speed of EDA-FRA* is much faster than that of the EDA-RA*. The pseudocode of the hybrid EDA-FRA* algorithm for fast minimum risk-cost path planning is shown in **Algorithm 2**.

The EDA algorithm independently process the cost data and output the best population of the feasible region (Line 5). The obtained best population is further processed to get the position of open points (Line 6-Line 9). The k-means method is then employed to obtain the centroid position of clusters, and to generate the heuristic distance factor and heuristic direction for RiskA*. The heuristic distance of RiskA* is improved from Euclidean distance to estimate the total risk cost from current node to destination, presented as

$$h_{\text{Dist}} = f_{heuDist}\sqrt{(x_D - x_i)^2 + (y_D - y_i)^2 + (z_D - z_i)^2} \tag{25}$$

$$f_{heuDist} = \min \left\{ \left( \frac{1}{N_{V_{open}}} \sum_{v_i \in V_{open}} C_{r\_v_i} \right), \left( \frac{1}{N_{V_{Ctrs}}} \sum_{v_j \in V_{Ctrs}} C_{r\_v_j} \right) \right\} \quad (26)$$

in which $h_{Dist}$ is the improved heuristic distance for RiskA*, and $(x_i, y_i, z_i)$ is the coordination of current point and $(x_D, y_D, z_D)$ is the destination point, which are used to compute the Euclidean distance. Note that $f_{heuDist}$ is the heuristic distance factor, which takes the minimum value among the mean risk cost of all open points $V_{open}$ and the mean risk cost of cluster centroid points $V_{Ctrs}$. Taking the minimum value is because smaller heuristic value leads to better quality of solution for RiskA* algorithm. The solution will reach optimum when heuristic value down to zero, and the RiskA* is then equivalent to Dijkstra. The determination of the heuristic value for the algorithm is to make a trade-off between the quality of solution and computational efficiency.

---

**Algorithm 2**: EDA-FRA* algorithm for fast minimum risk-cost 3D path planning

1: Result: *path*, *TotalCost*
2: load CostDataset
3: *path* = null
4: *TotalCost* = 0;
5: BestPop=EDA(CostDataset); % obtain the best population
6: **if** sum(BestPop(*i*, *j*, *k*))==1 % obtain open points for path searching
7:     IndivPosn=[*i*, *j*, *k*]; % get the index of the open points
8:     Posn=[Posn; IndivPosn];
9: **end**
10: [Ctrs, heuDist*f*, heuDrctn]=k-means(Posn); % obtain heuristic information
11: [*path*, *TotalCost*]=RiskA*(CostDataset, heuDist*f*, heuDrctn);
12: **Return** (*path*, *TotalCost*)

---

The total risk cost function from origin to destination can be then described as

$$f(C) = g(C) + h_{Dist}(C) = \int_{v_o}^{v_i} C(v) dv + h_{Dist}(C) \quad (27)$$

where $g(C)$ is the integral of the risk cost from origin $v_o$ to current point $v_i$, and $h_{Dist}(C)$ is obtained by Eq. (25). Note that $C(v)$ is the risk cost of path point $v$.

The heuristic direction is represented by a set of segments, which starts from origin connecting cluster centroids to the destination (Fig. 6(c)). The heuristic distance provides global information by selecting the node with the smallest distance $f(C)$ to destination. While the heuristic direction assists with local searching by evaluating the distance cost from current node to the local cluster centroid. The key part of RiskA* framework is given as **Algorithm 3**.

---

**Algorithm 3**: Improved RiskA* algorithm for minimum risk-cost path planning

1: Result: *path*
2: load CostDataset, heuDist*f*, heuDrctn
3: **for** *i* = 1:26 % possible choices of next position
4:     next = [*x*(*i*), *y*(*i*), *z*(*i*), CostDataset(*x*(*i*), *y*(*i*), *z*(*i*))]; % next position and its cost value
5:     Motion=[Motion; next];
6: **end**
7: MotionMode{} = Motion; % record the moving cost matrix for every open points
8: *g*(c)=*f*(MotionMode{c});
9: hDistDest(c)=*f*(heuDist*f*, *v*(c)); % heuristic distance to destination
10: hDistCtrs(c)=*f*(heuDrctn, *v*(c)); % heuristic distance along the cluster centriods to destination
11: **if** hDistCtrs(c)-hDistDest(c)<ε % degree of deviation from main track, ε>=0
12:     hDist(c)=hDistDest(c); % Dev is acceptable, put currrent point a small estimated cost
13: **else**
14:     hDist(c)=hDistCtrs(c); % unacceptable, put currrent point a high estimated cost, it will be then discarded
15: **end**
16: *f*(c)=*g*(c)+hDist(c); % distance funtion
17: **if** *f*(c) < open(*f*(:)) % distance from current point to destination less than that of the points in openlist
18:     Putting current point as father point, and it will be included in the path
19: **end**
20: **Return** (*path*)

The motion mode function is to obtain the risk cost value of the adjacent 26 possible positions. The cost matrix of all feasible points is computed and stored in the MotionMode array, which is employed to compute the total risk cost $g(C)$ from origin point to current point. For the second part, the heuristic distance $h_{\text{Dist}}$ is computed to estimate the distance from current point to destination. The $h_{\text{Dist}}^{Dest}$ (Line 9) is computed by Eq. (25). While $h_{\text{Dist}}^{Ctrs}$ (Line 10) the heuristic distance along the cluster centroids to destination is the sum of segment distance products the heuristic factor $f_{heuDist}$. As Fig. 6(c) depicted, the heuristic distance path is directly from current node to destination, while the heuristic direction path is generated by connecting the cluster centroids. The $h_{\text{Dist}}^{Dest}$ globally inspires the path searching process, and the $h_{\text{Dist}}^{Ctrs}$ assists with local searching by evaluating the deviation between heuristic direction and actual path searching one (Line 11 to Line 14). If the deviation exceeds the threshold $\varepsilon$, that node will be given an infinite cost and it will be removed in the next iteration. Here the $\varepsilon$ is taken as 0.2, which presents the deviation between main track and current track.

With the development above, the framework and relationship of EDA-RA* and EDA-FRA* are presented in Fig. 7. There are two similarities of the two hybrid algorithms. They use same risk cost data as input and they employ EDA to generate initial solution of feasible regions. The difference between the two algorithm is prominent. For EDA-RA*, the A* algorithm is called for each single iteration to produce path for all species. By giving dominant path points high probability of been selected into feasible region, the final optimized path will effectively search all space and obtain a good quality of solution. While for EDA-FRA*, the main loop is to optimize the feasible region. It calls RiskA* only once after obtaining the optimized feasible region. That significantly saves computational time whereas the optimized feasible region may not be as good as the EDA-RA* has.

The EDA-RA* can achieve better feasible region and solution quality, but it costs more time to compute. In comparison, the EDA-FRA* is much faster as it calls RiskA* only once while the solution quality may not compete with EDA-RA*. Simulation and case studies are conducted in next section to validate the proposed hybrid algorithms in term of computational efficiency and quality of solutions.

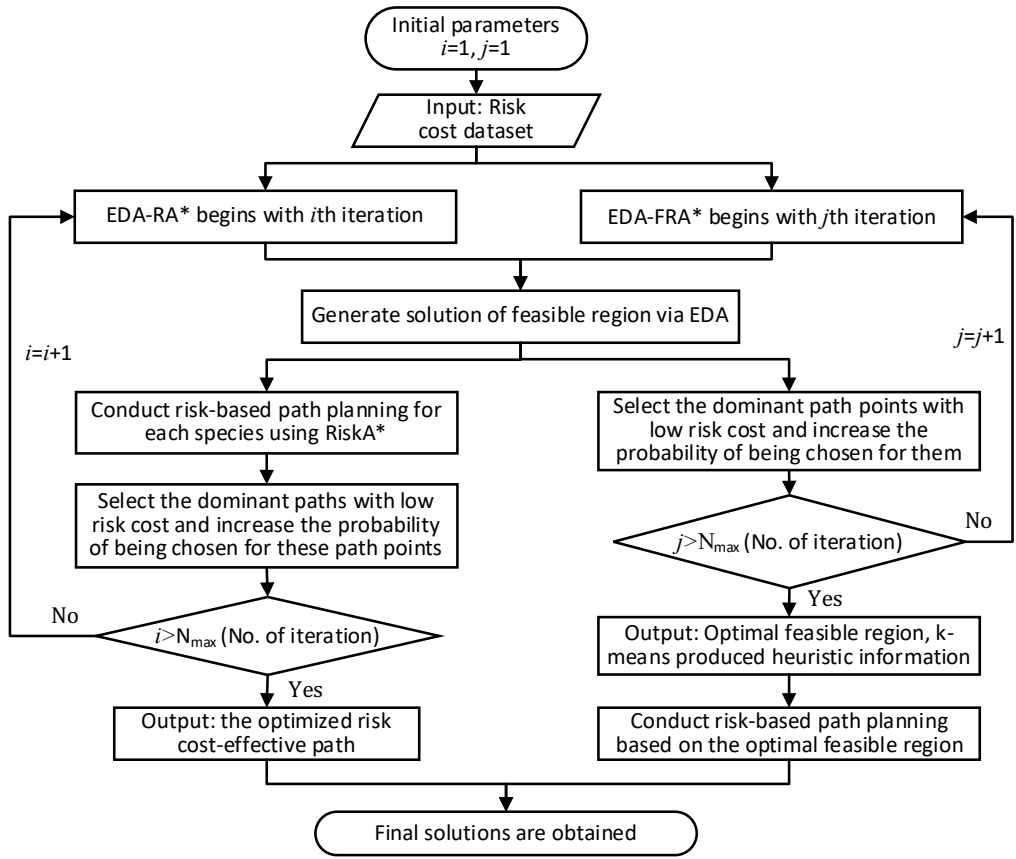

**Fig. 7**. Framework of EDA-RA* and EDA-FRA* algorithms.

## 5. Simulation studies

To validate the proposed risk cost assessment model and the developed hybrid risk-based 3D path planning algorithms, we perform simulations and case studies in the context of a representative metropolitan area. Firstly, the risk assessment model is implemented in a real-world environment to generate 3D risk-aware airspace map. We then apply the proposed hybrid algorithms to the generated 3D airspace to produce the risk-cost-effective path. Lastly, we conduct simulations and statistical analysis to test how well the proposed risk assessment model and algorithms can be generalized to other urban patterns.

### 5.1. Case study of risk cost assessment model

A typical metropolitan area (6km×6km) in Singapore is selected for the modelling of risk aware airspace, and the allowable altitude in this study is chosen as 120 meters (400 feet) above the ground. The size of each air block is 100m×100m×30m. The selected metropolitan area has dense high-rise buildings, shopping centers, city squares, residential areas with dense population, parks, etc., which are representative for modern mega cities. The selected environment has two administrative districts, and the average population densities are $8.358 \times 10^3$ and $7.219 \times 10^3$ (people/km$^2$) (WorldoMeter 2021). The average traffic density in the given area is obtained as $7.12 \times 10^3$ (vehicle/km$^2$) (SG Land Transport Authority 2021). Based on the average population and vehicle densities, we can estimate the fine population density distribution using Eq. (10) and traffic density distribution using Eq. (11). Obtained population and traffic distribution results are demonstrated in Fig. 8.

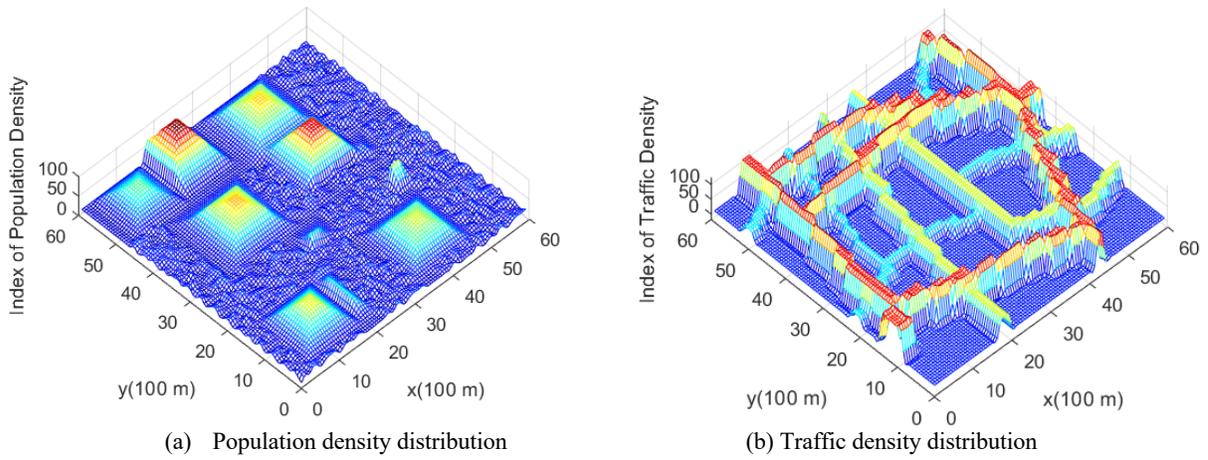

(a) Population density distribution  (b) Traffic density distribution

**Fig. 8.** Distribution of population density and traffic density in selected urban environments.

In this case study, the UAV is selected as one of the most commonly used drones (DJI Phantom 4). The weight is 1.38kg, and the crash probability $P_{\text{crash}}$ is $6.04 \times 10^{-5}$ per flight hour (Shaokun 2018). The explored area of UAV impacts the ground is $S_{\text{hit}} = 0.0188$ m$^2$ and the drag coefficient is $R_I = 0.3$ (Koh et al. 2018). The number of casualties caused by average traffic accident is $R_f^v = 0.27$ (Budget Direct Insurance 2021). For the integrated risk cost model, the weight factors of fatality risk, property damage and societal impact are given as $\alpha_1 = 0.5, \alpha_2 = 0.25$ and $\alpha_3 = 0.25$, respectively (see Eq. (19)). The fatality risk cost is given a high weight with 50% of the total weight, while the property damage cost and societal impact cost are given the same weight with 25% each in this study. Based on the obtained data, the total integrated risk cost of each flight layer is computed and demonstrated in Fig. 9.

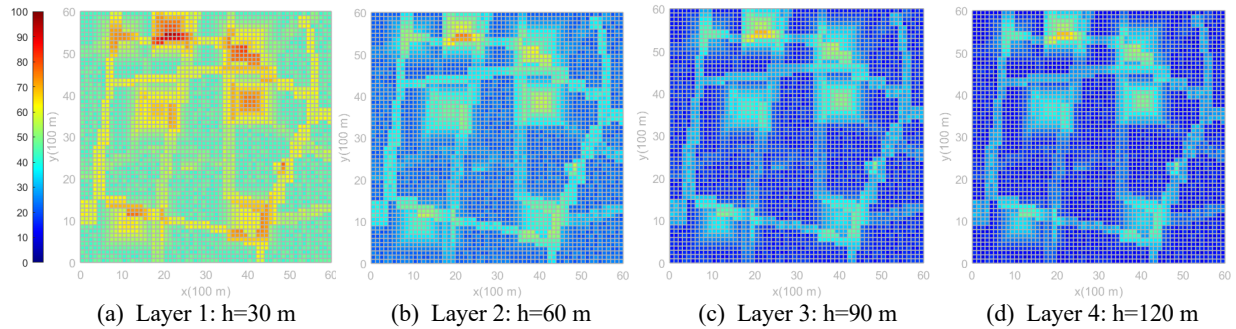

(a) Layer 1: h=30 m  (b) Layer 2: h=60 m  (c) Layer 3: h=90 m  (d) Layer 4: h=120 m

**Fig. 9.** Risk cost mapping for different flight altitude.

In Fig. 9(a), the altitude of the flight layer is 30 meters and the average risk cost of whole area is high. The areas with high fatality risks are well identified in the map. For instance, the location (20, 55) on the map has the highest risk cost, as there are shopping streets, highway intersections, and dense population in real-world environments. Thus, the fatality risk and property damage risk in there are high, making the total risk cost high. In Layer 2 (Fig. 9(b)), with the increase of flight altitude, risk costs are significantly reduced for property damage and noise impact as denoted in Eq. (14) and Eq. (17), whereas the fatality risk cost increases 7.7% compared with the one in the Layer 1.

In the third and fourth flight layers as demonstrated in Fig. 9(c) and Fig. 9(d), the flight altitude increases to 90 and 120 meters. The high-risk areas in these two layers are still clearly identified while the total risk cost has not changed much from Layer 3 to Layer 4 because of two reasons. For one thing, the fatality risk only slightly increases (4.09% from Layer 2 to Layer 3 and 2.66% from Layer 3 to Layer 4) after the altitude passing 60 meters, as the impact damage over such height is mostly the same, which is causing fatalities. For another, the influence of the societal impact exceeds the height threshold (40 meters in Eq. (17) and Fig. 5(b)) and contributes nothing to the risk cost, while the property damage cost is significantly small. The risk cost in Layer 3 and Layer 4 are therefore significantly small while with high risk areas being clearly identified.

*5.2. Risk-based path planning analysis*

With the risk-aware 3D airspace map, the risk-based path planning is conducted using EDA-RA* and EDA-FRA* algorithms. For comparison, Dijkstra and ACO algorithms are employed in the same environments. What is more, we investigate the impact of different risk types on the risk-based path planning and safe UAV operation.

A. *3D risk-based path planning in real-world environment*

In general, the produced 3D risk-cost-effective paths are able to avoid obstacles and high-risk areas identified by our proposed risk assessment model. The results of 3D view with risk map, without risk map and top view are presented in Fig. 10(a), Fig. 10(b) and Fig. 10(c). Observed from 3D view, the drone flight height for most of the time is 120 meters, which is the top layer of the modelled environments. Flying at such height significantly reduced the property damage risk and noise impact cost, and the fatality cost can also be reduced by avoiding high population density and vehicle density areas such as locations (10, 12), (40, 40) and (42, 39) shown in the top view.

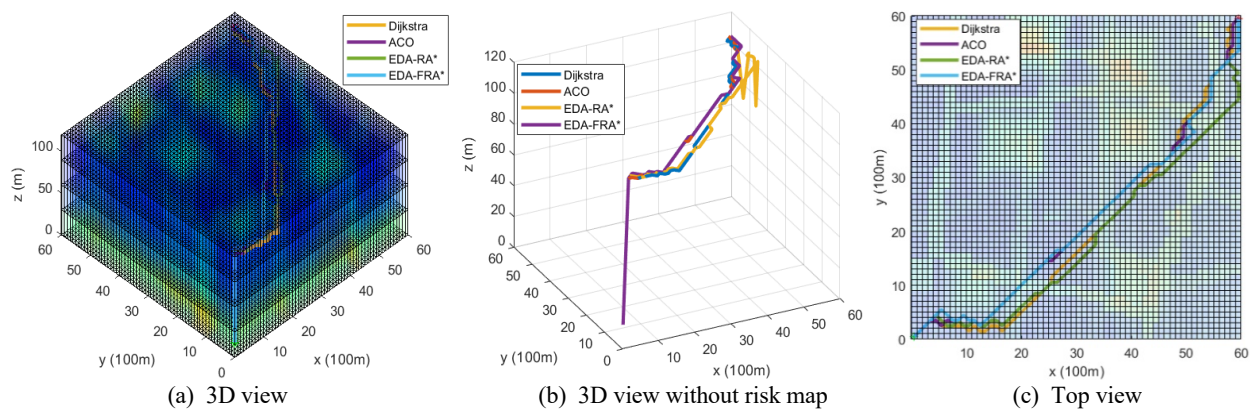

(a) 3D view  (b) 3D view without risk map  (c) Top view

**Fig. 10.** Results of 3D risk-based path planning in real-world environments.

In real-world environment, these identified high risk locations are shopping street, hospital, school or highway conjunctions, where the population density and vehicle density are significantly higher than the rest of the areas. Being able to quantitively identify high risk areas using our proposed model facilitates the risk management of low altitude urban airspace and risk-based path planning. Which subsequently enables safe UAV operations in metropolitan environments.

The path planning results of four algorithms are obtained and presented in Fig. 11 and

Table 1. Compared with Dijkstra algorithm, EDA-FRA* produced path has 2.06% shorter distance while using a mere of 3.05% computational time, whereas the risk cost of the path is 4.58% greater than the Dijkstra one. Followed by the path produced by ACO with 1.51% more cost and 1.22% longer distance. Its computational time is, however, dramatically greater than all of other three algorithms, with 785.98% greater than Dijkstra method. For EDA-RA* algorithm, its performance makes a good trade-off between risk cost, path distance and computational time.

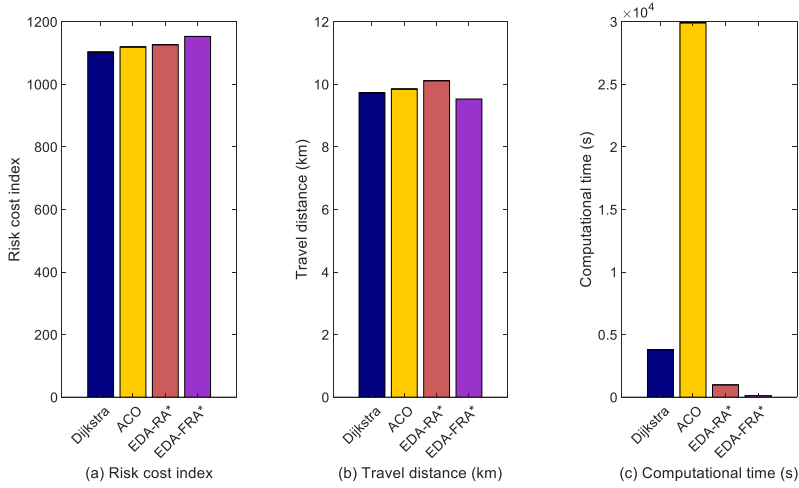

**Fig. 11.** Comparison results of four path planning algorithms.

**Table 1**
Detailed performance results of the four path planning algorithms.

| Performance\Indicators | Computational results | | | | Performance ratio compared with Dijkstra method | | |
|---|---|---|---|---|---|---|---|
| | Dijkstra | ACO | EDA-RA* | EDA-FRA* | ACO | EDA-RA* | EDA-FRA* |
| Risk cost index | 1102.60 | 1119.32 | 1126.26 | 1153.07 | 101.51% | 102.14% | 104.58% |
| Distance (km) | 9.73 | 9.85 | 10.11 | 9.53 | 101.22% | 103.94% | 97.94% |
| Computation time (s) | 3806.21 | 29916.15 | 984.87 | 116.22 | 785.98% | 25.88% | 3.05% |

*B. Impact of different risk types on risk-based path planning*

In order to demonstrate the impact of different risk types on safe UAV operations, we conduct four groups of path planning simulation with risk-related variables controlled and all constrains applied (see Section 4.1). Four simulations are planned as: (1) Path1: without consider any risk; (2) Path2: only consider fatality risk; (3) Path3: consider fatality risk and property damage risk; (4) Path4: consider all three risks. The environment of the four simulations is the same, which is generated in Section 5.1. The EDA-FRA* algorithm performs well among the benchmark methods in terms of computational efficiency and effectiveness and is therefore applied for each of the four simulations. Obtained results are shown in Fig. 12 and Table 2.

The Path1 goes from origin (1, 1, 1) to destination (60, 60, 4) with almost a straight line. This path avoids obstacles in positions like (30, 32) and (40, 40) shown in Fig. 12. However, it does not avoid high population density and vehicle density areas where the risk costs are high, resulting in the total risk cost of Path1 is the highest among all paths, with risk cost index of 1698.25 (see Table 2). Whereas the distance of the path is the shortest as it goes almost straightly. For Path2 the fatality risk is taken into account. This path successfully avoids the high population density areas (10, 12), (40, 45), (42, 15) shown in Fig. 8(a) and high vehicle density areas (8, 20), (45, 10) and (55, 40) shown in Fig. 8(b). As the fatality risk cost has a high proportion in the total risk cost model, avoiding high fatality risk areas makes a significant reduction (23.68%) of risk cost for Path2, compared with Path1. While the distance of Path2 is the greatest among the four paths, as more distance is travelled to avoid high risk areas. For Path3 the fatality risk and property damage risk are both considered. The produced path3 not only averts the high fatality risk areas, but avoid dense high-rise building areas like (20, 20) and (30, 30) shown in Fig. 4(a), which Path1 and Path2 fail to do so (see Fig. 12). By adding property damage risk into the model, the Path3 is able to additionally avoid dense building areas, thus the risk cost of its path further drops by 5.94% compared with Path2. For Path4, it takes all three risk types into consideration. However, the most part of the path is at flight level 120 meters where the property damage and noise impact have tiny contributions to the risk cost. So, the risk cost of Path4 has no remarkable reduction (a mere of 2.48%) compared with Path3.

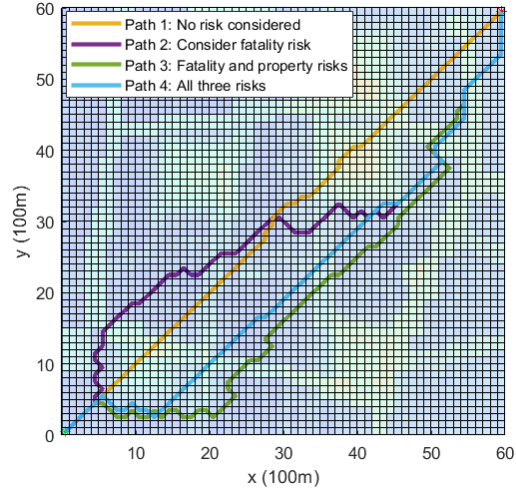

**Fig. 12.** Top view of produced paths for considering four different risk types.

Overall, the fatality risk is identified as the primary factor in the risk assessment model and in the subsequent cost-effective path planning, followed by property damage risk and noise impact cost. With more risk types getting considered, the more accurate risk aware 3D airspace map will be generated to produce lower risk cost path for safe UAV operations. While the distance of the produced path will be greater as the UAV needs more movements to avoid obstacles and high-risk areas.

**Table 2**
Results of the four produced paths considering different risk types.

| Performance\Indicators | Computational results | | | | Performance ratio compared with Path1 | | | |
|---|---|---|---|---|---|---|---|---|
| | Path1 | Path2 | Path3 | Path4 | Path1 | Path2 | Path3 | Path4 |
| Risk cost index | 1698.25 | 1296.02 | 1195.18 | 1153.10 | 100.00% | 76.32% | 70.38% | 67.90% |
| Travel Distance (km) | 8.35 | 10.22 | 10.01 | 9.53 | 100.00% | 122.32% | 119.88% | 114.08% |

*5.3. External validity of the risk assessment model*

To validate how well the proposed risk assessment model can be generalized to other urban patterns in mitigating third-party risk, we conduct external validity and randomly generate the parameters of 100 different urban patterns. We take population density from the integer range of $[5, 25] \times 10^3$ (people/km$^2$), which covers the most densely populated cities worldwide (Wikipedia 2021). The average traffic density is given as same as above in Section 5.1*5.1*. The building height distribution of all generated patterns follows log-normal distribution. The scope of the validation environment is 6km×6km×120m with size of each unit is 100m×100m×30m. Using the risk assessment model in Section 3 with the generated parameters, we obtain the risk cost value in each airspace unit for 100 independent simulation environments.

Each simulation has been independently conducted standard Dijkstra and risk-based Dijkstra via MATLAB software on a desktop equipped with an InteI E5-2680 @2.4Ghz CPU. The standard Dijkstra algorithm is performed in a normal map without considering third-party risk, while risk-based Dijkstra is conducted based on the risk map generated by our assessment model. Comparison is made between these two ways to see how much percentage of risk being mitigated by using risk-based method. The simulation starts from origin (1, 1, 1) to destination (60, 60, 4) to produce the risk-cost-effective path. Total risk cost for each simulation is obtained and presented in Fig. 13.

In the 100 generated samples (urban patterns), the mitigated risk of flight path in each individual urban pattern is greatly less than the unmitigated one produced by standard Dijkstra. To test how effective the risk can be mitigated for the population (all types of urban patterns), we conduct statistical analysis to find a 95% confidence interval for the percentage of the risk being mitigated by.

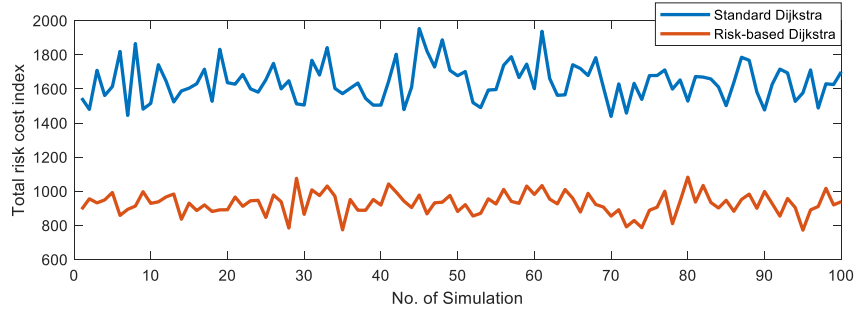

**Fig. 13.** Risk being mitigated by using risk-based method.

In this case, we have two sample groups and they are mitigated risk group (Group 1) and unmitigated risk group (Group 2). Each group has 100 samples (total risk cost index). As the sample sizes are large ($n_1 > 30$ and $n_2 \gg 30$), we use normal distribution the compute the confidence interval. The sample means ($\bar{x}_1$ and $\bar{x}_2$) and sample variances ($s_1^2$ and $s_2^2$) of the two groups are computed and presented in Table 3. The population means of Group 1 and Group 2 are presented as $\mu_1$ and $\mu_2$. The confidence interval of risk mitigation effect can be then described as $(\mu_2 - \mu_1)/\bar{x}_2$, where the interval of $\mu_2 - \mu_1$ can be computed as $(\bar{x}_2 - \bar{x}_1) \pm Z_{\alpha/2}\sqrt{s_1^2/n_1 + s_2^2/n_2}$. The obtained result shows that a 95% confidence interval for risk mitigation effect is $(\mu_2 - \mu_1)/\bar{x}_2 \in [0.4264, 0.4415]$. Which means that our proposed risk assessment model is effective for all types of urban patterns, and the average total risk can be mitigated by [42.64%, 44.15%] at 95% confidence level.

**Table 3**
Statistical analysis parameters of risk mitigation effect.

| Groups | Mitigated risk (Group 1) | Unmitigated risk (Group 2) |
| --- | --- | --- |
| Sample size | $n_1 = 100$ | $n_2 = 100$ |
| Sample mean | $\bar{x}_1 = 18584$ | $\bar{x}_2 = 32831$ |
| Sample variance | $s_1^2 = 1594657$ | $s_2^2 = 4859067$ |

*5.4. Reliability validation of the proposed algorithms*

Above we validated that the proposed risk assessment model is effective to capture risk features and mitigate total risks in all types of urban patterns. To test the reliability of proposed risk-based algorithms in solving path planning problems with different urban patterns and risk maps, we perform the simulations based on the urban patterns and risk cost value generated in the Section 5.3. the comparison involves three algorithms of risk-based Dijkstra, EDA-RA* and EDA-FRA*. Total risk cost and total distance of the produced path, as well as computational time are obtained and presented in Fig. 14 and Table 4.

In the overall 100 simulations, the Dijkstra algorithm performs the best in terms of the total risk cost of the produced path, followed by EDA-RA* and EDA-FRA* with average performance rate of 102.76% and 105.47%, respectively. It means that the average total risk cost of the path produced by EDA-RA* and EDA-FRA* algorithms are 2.76% and 5.47% greater than that of the Dijkstra method. In contrast, the EDA-FRA* algorithm provides the best performance in terms of average total distance and average computational time, with 0.4% shorter in total distance and a mere of 2.07% in computational time. Followed by EDA-RA* algorithm, it saves 0.17% travel distance while spends 24.12% computational time compared with Dijkstra method. Another interesting trend can be observed is that total risk cost (Fig. 14 (a)) and total travel distance (Fig. 14 (b)) are two variables with opposite trend. One variable increases resulting in decrease of the other (e.g. the 7th, 13th and 37th simulations). Which implies that to achieve a low risk cost path in risk-based environment, UAV always has to travel more distance to avoid high risk areas.

Deviation analysis of the 100 simulations with randomly generated risk environments shows that the proposed EDA-RA* and EDA-FRA* algorithms has good robustness, with a mere of 1.29% and 2.54% standard deviation in computing total risk cost, while 2.18% and 4.29% for total distance (Table 4). The results present that the proposed algorithms are reliable to solve complex risk-based 3D path planning problems for different environments.

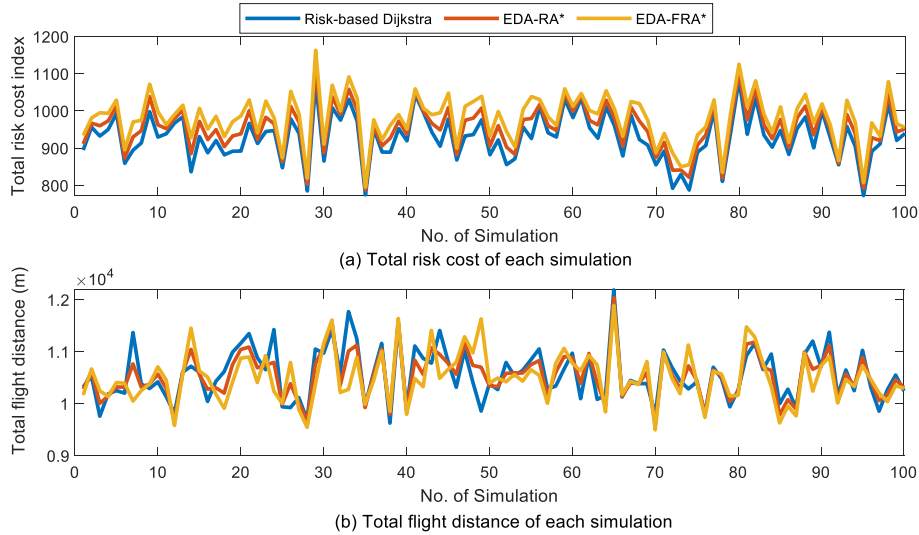

**Fig. 14.** Reliability testing results for proposed algorithms in terms of total risk cost and total flight distance.

In summary, the EDA-FRA* is the fastest algorithm in terms of computation and it is much faster than the other two algorithms. While it also has a good performance in obtaining the quality of solution. Its performance is made possible by the heuristic information scheme and the algorithm's structure (see Algorithm 3) without need to call A* function in every loop. The EDA-FRA* algorithm can be therefore used for more time-sensitive risk-based path planning missions. The EDA-RA* algorithm performs well in both quality of solution and computational efficiency, so it can be applied to the trade-off case. As an exact method, Dijkstra can be employed to produce the minimum risk cost path if computational time is not sensitive.

**Table 4**
Reliability testing results of the proposed algorithms with 100 independent simulations.

| Performance\Algorithms | Average total risk cost | | Average total distance (km) | | Average computational time (s) | Average percentage of the time |
| --- | --- | --- | --- | --- | --- | --- |
| | Indicator ratio | Standard Deviation | Indicator ratio | Standard Deviation | | |
| Dijkstra | 100.00% | 0.00% | 100.00% | 0.00% | 4120.35 | 100.00% |
| EDA-RA* | 102.76% | 1.29% | 99.83% | 2.18% | 993.68 | 24.12% |
| EDA-FRA* | 105.47% | 2.54% | 99.60% | 4.29% | 85.24 | 2.07% |

## 6. Conclusions

In this article, we investigated the third party risk assessment model and risk-based path planning problems for safe UAV operations in metropolitan environments. Main findings are concluded as follows.

(1) The proposed integrated risk assessment model is able to capture comprehensive risk types in urban environments including fatality risk, property damage risk and societal impact risk. The statistical analysis results show that the proposed model is effective for all types of urban patterns. The model enables quantitative risk assessment for urban airspace, which facilitates the risk aware airspace modelling and risk-based path planning.

(2) The introduced gravity model for population density and vehicle density estimation can identify high population density areas in a finer scale. Which means that the population density in one district will not be averagely treated but being identified with actual high-density areas while releasing low density areas for UAV operations.

(3) The developed hybrid 3D risk-based path planning algorithms outperform the existing methods in terms of solution quality and computational efficiency. The proposed hybrid EDA-FRA* algorithm performs best in computational time while it can still reach an average of 94% optimality of solution quality. Besides, the proposed EDA-RA* algorithm has better performance in making trade-off between solution quality and computational time.

(4) With more risk sources been considered, the total risk cost of produced path deceases. That is because the risk-based path planning algorithm will avoid areas with high risks, provided the risk types are considered and assessed. That further justifies the significance of our proposed integrated risk assessment model, which is able to incorporate more risks with different types.

The work of this article can be further improved from several aspects. For instance, a more accurate estimation model for population density could be developed if multisource data is available (Deville et al. 2014). Besides, collision risk between UAV and manned aircraft can also be incorporated into the integrated model for these environments where manned and unmanned aircraft are allowed to co-exist.

**Acknowledgement**

The authors would like to thank Dr. Yu Wu and Dr. C.H. John Wang for their insightful suggestions on risk modelling and algorithm design. This research is supported by the Civil Aviation Authority of Singapore and the Nanyang Technological University, Singapore under their collaboration in the Air Traffic Management Research Institute. Any opinions, findings and conclusions or recommendations expressed in this material are those of the authors and do not reflect the views of the Civil Aviation Authority of Singapore.